\documentclass[12pt]{article}

\usepackage{amssymb,epsfig,amsmath,amsthm}

\newcommand{\nc}{\newcommand}    \newcommand{\rnc}{\renewcommand}
\nc{\nf}{\newfont}
\nc{\al}{\alpha}                 \nc{\be}{\beta} \nc{\ga}{\gamma}                 \nc{\Ga}{\Gamma} \nc{\De}{\Delta}
\nc{\de}{\delta} \nc{\ep}{\epsilon}               \nc{\ze}{\zeta} \nc{\ka}{\kappa}                 \nc{\la}{\lambda}
\nc{\om}{\omega}                 \nc{\Om}{\Omega} \nc{\si}{\sigma}                 \nc{\tht}{\theta}
\nc{\vep}{\varepsilon} \nc{\La}{\Lambda}
\usepackage{dsfont}
\nc{\ds}{\mathds}

\nc{\nin}{\noindent}             \nc{\nn}{\nonumber}
\nc{\ul}{\underline}

\nc{\Ra}{\Rightarrow} \nc{\ra}{\rightarrow} \nc{\lra}{\longrightarrow}

\nc{\sfrac}[2]{\mbox{\small{$\frac{#1}{#2}$}\normalsize}}
\nc{\vc}[1]{\stackrel{#1}{\small{\sim}}}
\nc{\add}{\qq\mb{and}\qq}

\nc{\chapterbib}{\rnc{\bibname}{References} \bibliography{../extras/phdrefs} \bibliographystyle{acm}}  
\nc{\bcen}{\begin{center}}       \nc{\ncen}{\end{center}}
\nc{\bfig}{\begin{figure}}       \nc{\nfig}{\end{figure}}
\nc{\beqn}{\begin{equation}}     \nc{\neqn}{\end{equation}}
\nc{\beqns}{\begin{equation*}}  \nc{\neqns}{\end{equation*}}
\nc{\beqnays}{\begin{eqnarray*}} \nc{\neqnays}{\end{eqnarray*}}
\nc{\beqnay}{\begin{eqnarray}}   \nc{\neqnay}{\end{eqnarray}}
\nc{\barr}{\begin{array}}        \nc{\narr}{\end{array}}
\nc{\bcoll}{\begin{collarray}}   \nc{\ncoll}{\end{collarray}}
\nc{\bcolls}{\begin{collarray*}} \nc{\ncolls}{\end{collarray*}}
\nc{\btab}{\begin{tabular}}      \nc{\ntab}{\end{tabular}}
\nc{\bdes}{\begin{description}}  \nc{\ndes}{\end{description}}
\nc{\ben}{\begin{enumerate}}     \nc{\nen}{\end{enumerate}}
\nc{\biz}{\begin{itemize}}       \nc{\niz}{\end{itemize}}
\nc{\bqt}{\begin{quote}}     \nc{\nqt}{\end{quote}}
\nc{\bpic}{\begin{picture}}  \nc{\npic}{\end{picture}}
\nc{\mkbx}{\makebox}

\nc{\bmat}{\begin{pmatrix}} \nc{\nmat}{\end{pmatrix}}
\nc{\bals}{\begin{align*}}
\nc{\nals}{\end{align*}}
\nc{\bcase}{\begin{cases}}  \nc{\ncase}{\end{cases}}

\newenvironment{arabiclist}{%

    \let\item\Item
    \begin{enumerate}
    }{%
    \end{enumerate}}

\nc{\bal}{\begin{arabiclist}}
\nc{\nal}{\end{arabiclist}}

\nc{\ch}{\chapter}
\nc{\sn}{\section}
\nc{\ssn}{\subsection}           \nc{\sssn}{\subsubsection}
\nc{\sep}{\bcen ---oOo--- \ncen}
\nc{\lrb}{\left (}      \nc{\rrb}{\right )}     
\nc{\lcb}{\left \{}     \nc{\rcb}{\right \}}    
\nc{\lsb}{\left [}      \nc{\rsb}{\right ]}     
\nc{\lmb}{\left |}    \nc{\rmb}{\right |}   
\nc{\lfb}{\left \lfloor}      \nc{\rfb}{\right \rfloor}     
\nc{\lceilb}{\left \lceil}      \nc{\rceilb}{\right \rceil}     
\nc{\blcb}{\big \{} \nc{\brcb}{\big \}} \nc{\blrb}{\big (} \nc{\brrb}{\big )}
\nc{\rb}[1]{\raisebox{1.5ex}{#1}}
\nc{\hf}{\hfil}                  \nc{\hfl}{\hfill}
\nc{\fn}{\footnote}              \nc{\ph}{\phantom}
\nc{\fb}{\fbox}
\nc{\hb}{\hbox}                  
\nc{\vb}{\vbox}                  \nc{\vt}{\vtop}
\nc{\q}{\quad}                   \nc{\qq}{\qquad}
\nc{\hs}{\hspace}                \nc{\vs}{\vspace}
\nc{\hsst}{\hspace*}             \nc{\vsst}{\vspace*}
\nc{\ssk}{\smallskip}            \nc{\msk}{\medskip}
\nc{\bsk}{\bigskip}
\nc{\nl}{\newline}               \nc{\np}{\newpage}
           \nc{\corr}{\mbox{corr}}

\nc{\scsz}{\scriptsize}      \nc{\fnsz}{\footnotesize}
\nc{\tm}{truemm}                 \nc{\tp}{truept}
\nc{\half}{\sfrac{1}{2}}

\nc{\clearemptydoublepage}{\newpage{\pagestyle{empty}\cleardoublepage}}
\nc{\R}{\mathbb{R}} \nc{\Z}{\mathbb{Z}}  \nc{\C}{\mathbb{C}} \nc{\N}{\mathbb{N}}

\nc{\Var}{\mathrm{Var}}
\nc{\Exp}{\mathbb{E}} \nc{\Cov}{\mathrm{Cov}} \nc{\Prb}{\mathbb{P}} \nc{\tri}{\mb{tri}} \nc{\rect}{\mb{rect}}
\nc{\step}{\mb{step}} \nc{\sinc}{\mb{sinc}} \nc{\vect}{\mb{vec}} \nc{\Po}{\mathrm{Po}} \nc{\Bi}{\mathrm{Bi}}
\nc{\mult}{\mathrm{Mult}} \nc{\Un}{\mathrm{U}} \nc{\wt}{\widetilde}
\makeatletter
\def\imod#1{\allowbreak\mkern10mu({\operator@font mod}\,\,#1)}
\makeatother

\nc{\sgn}{\mathrm{sgn}}
\nc{\lcm}{\mathrm{lcm}}

\nc{\abs}[1]{\lvert #1 \rvert} \nc{\norm}[1]{\lVert #1 \rVert}

\nc{\mc}{\mathcal} \nc{\itxt}{\intertext}  \nc{\E}{\mathbb{E}}  \nc{\eqd}{\stackrel{\mc D}{=}}
\nc{\ebd}{:=} \nc{\gl}{\gtrless} \nc{\ind}{\mathbb{I}}

\nc{\PBD}{\mathrm{PBD}} \nc{\B}{\mathrm{Bi}} \nc{\PB}{\mathrm{PB}}

\nc{\pd}{\partial}
\nc{\bs}{\ensuremath{\boldsymbol}}
\nc{\mb}{\ensuremath{\mathbf}}

\rnc{\eqd}{\stackrel{\mc D}{=}}
\nc{\Al}{\bs \alpha} \nc{\ainf}{\Alpha_{\infty}} \nc{\td}{\stackrel{\mc D}{\lra}} \nc{\Oma}{\bs{\om}}
\nc{\oinf}{\Oma_{\infty}} \rnc{\tp}{\stackrel{\Prb}{\lra}}

\nc{\etal}{{\it et al }} \nc{\eg}{{\it e.g. }} \nc{\egns}{{\it e.g.}} \nc{\ie}{{\it i.e. }} \nc{\ci}{($i$)}
\nc{\cii}{($ii$)} \nc{\ciii}{($iii$)} \nc{\etc}{{\it etc. }} \nc{\aka}{{\it a.k.a. }} \nc{\apriori}{{\it a priori }}
\nc{\ou}{Ornstein-Uhlenbeck }

 \nc{\X}{\bs X} \nc{\x}{\bs x} \nc{\pdv}{\partial} \nc{\dif}{\mathrm{d}}
\nc{\y}{y} \nc{\W}{W} \nc{\I}{\ensuremath{\mathrm{I}}} \nc{\M}{\bs M}

\newtheorem{thm}{Theorem}[section]
\newtheorem*{thm*}{Theorem}
\newtheorem{prop}[thm]{Proposition}
\newtheorem{lem}{Lemma}[section]

\theoremstyle{definition}

\newtheorem{rem}{Remark}[section]

\newtheorem*{example*}{Example}

\numberwithin{equation}{section}

\makeatletter
\def\section{\@startsection{section}{1}{\z@}%
                           {-3.5ex \@plus -1ex \@minus -.2ex}%
                           {2.3ex \@plus.2ex}%
                           {\Large\bfseries\centering}}

\def\subsection{\@startsection{subsection}{1}{\z@}%
                           {-3.5ex \@plus -1ex \@minus -.2ex}%
                           {2.3ex \@plus.2ex}%
                           {\large}}

\title{Limit theorems for projections of random walk on a hypersphere}
\author{Max Skipper}

\begin{document}

\maketitle


\begin{abstract}
We show that almost any one-dimensional projection of a suitably scaled random walk on a hypercube, inscribed in a
hypersphere, converges weakly to an Ornstein-Uhlenbeck process as the dimension of the sphere tends to infinity. We
also observe that the same result holds when the random walk is replaced with spherical Brownian motion. This latter
result can be viewed as a ``functional'' generalisation of Poincar\'e's observation for projections of uniform measure
on high dimensional spheres; the former result is an analogous generalisation of the Bernoulli-Laplace central limit
theorem. Given the relation of these two classic results to the central limit theorem for convex bodies, the modest
results provided here would appear to motivate a functional generalisation.
\end{abstract}

\section{Introduction}

\nc{\sd}{^{(d)}}
\nc{\Sd}{\mc S^{d-1}_d} \nc{\So}{\mc S^{d-1}_1} \nc{\bth}{\bs \tht}

Let $\mc S^{d-1}_r$  be the spherical surface, centered at the origin, of radius $\sqrt r$ and let $\X$ be uniformly
distributed on $\Sd$. A classic observation dating back to Maxwell, Poincar\'e and Borel is that
the distributions of the first $K$ coordinates of $\X$ converge to independent standard normals as $d \to \infty$; see
Diaconis and Freedman \cite{df87} for an historical account. Recently there has been much work done on a generalisation
of this result which seeks to replace $\Sd$ with an arbitrary convex body $\mc K \in \R^d$ and
the coordinates of $\X$ with arbitrary linear projections $\langle \bs \tht, \X \rangle$, $\bth \in \So$. If $\X$ is
now uniformly distributed on $\mc K$, what is now known as the central limit theorem for convex bodies asserts that,
under suitable conditions on $\mc K$, the law of $\langle \bs \tht, \X \rangle$ is {\it approximately} Gaussian for
{\it most} $\bth \in \So$. We refer to Klartag \cite{klartag07} for precise statements of the theorem complete with
quantitative definitions of the words ``approximately'' and ``most'', as well as an overview of previous work; see
Milman \cite{milman09} for a more recent account including improved estimates for special cases.

Relative to the central limit theorem for convex bodies, the results presented in this note take the classical
Maxwell-Poincar\'e-Borel observation in another direction --- replacing an observation of projections of uniform
measure on $\Sd$ with an observation of projections of two particular ``uniform'' processes on $\Sd$. The two processes
considered are spherical Brownian motion (SBM) on $\Sd$ and a nearest neighbour random walk on the hypercube
$\{-1,1\}^d$ (inscribed in $\Sd$); they are referred to as ``uniform'' only because their invariant measures are
uniform on their support. Assuming $\bth\sd \in \So$ and letting $\X\sd$ denote either of the above mentioned processes
started at $\x\sd$, our main result states that if $|\bth\sd|_\infty \to 0$ and $\langle \bth\sd,\x\sd\rangle \to u$,
then $\langle \bth\sd,\X\sd\rangle$ converges weakly to an Ornstein-Uhlenbeck (OU) process $U$ started at $u$. (The
condition on the $\infty$-norm of $\bth\sd$ is unnecessary in the SBM case.) Just as the Maxwell-Poincar\'e-Borel
observation represents a special case of a ``na\"{\i}ve'' (non-quantitative) central limit theorem for convex bodies,
this modest result represents a special case of a na\"{\i}ve {\it functional} central limit theorem for convex bodies
--- developments and applications of which we hope to report in a subsequent paper.

Indeed, the original impetus for this work was derived from the practical problem of how to extract macroscopic
dynamics from a randomly evolving system where an explicit microscopic description is given. Typically, the microscopic
behaviour is modelled by a large system of coupled stochastic differential or difference equations driven by continuous
or discrete Markov processes. In contrast, the dynamics of interest are those of a smaller number of functionals of the
microscopic variables which are, in general, non-Markov. Since usually solutions must be obtained numerically, the main
objective is to find a self-contained approximate description of the sought-after dynamics without needing to fully
resolve the dynamics of the larger system; see Givon {\it et al} \cite{givon} for an informative survey.

A particularly relevant example is the Ehrenfest model of heat exchange between two isolated bodies, first published in
1907 in an effort to reconcile the irreversibility and recurrence in Boltzmann's kinetic theory of gases. (See Tak\'acs
\cite{takacs79} for an historical account of early work and Kac \cite{kac} for a discussion of Zermelo's
irreversibility/recurrence paradox.) The original model involves $d$ balls --- representing energized gas molecules ---
distributed among two urns --- the isolated bodies. The microscopic dynamics are such that at each time increment a
ball is drawn out at random and placed in the opposite urn from whence it came. The macroscopic variable of interest is
the number of balls in the first urn.

As is well-known, one may describe the allocation of the balls in the Ehrenfest model by the vector $\X(n) \in
\{0,1\}^d$ where $X_i(n) = 1$ if the $i$th ball is in the first urn after $n$ transitions and $X_i(n) = 0$ otherwise.
Moreover, the Ehrenfest dynamics imply that $\X = \{\X(n)\}_{n \in \N_0}$ is a random walk on the hypercube
$\{0,1\}^d$. What is special about this example is that because the $d$ microscopic variables (balls) are exchangeable,
the macroscopic Ehrenfest process $\sum_{i=1}^d X_i$ is also Markov and hence an exact self-contained description is
readily obtained. Nevertheless, it wasn't until $40$ years after the publication of the model that Kac \cite{kac}
managed to derive the transition probabilities. As part of his work Kac found the transition probabilities of a
suitably normalised Ehrenfest process and provided a sketch of how they converge to those of an OU-process as
$d\ra\infty$.

As far as we know, the most general extension of the Ehrenfest model that has some overlap with the work here is that
given by Schach \cite{schach}. Schach's model consists of $d$ balls distributed among $K$ urns where at each transition
a ball is moved from urn $j$ to urn $k$ with probability proportional to the number of balls in urn $j$ and a given
number $p_{jk}$. Again, the microscopic variables are (the locations of) the balls and the macroscopic variables are
the numbers of balls in each urn. Again, the macroscopic variables are Markov. Schach shows that as $d\ra \infty$ the
suitably normalised $K$-variate macroscopic process converges weakly to a $K$-variate OU-process. He also includes an
account of earlier work and discusses applications of his results.

The Ehrenfest models are examples of models in which the (normalised) macroscopic process retains the Markov property
and may be reasonably approximated by a diffusion process. Since the publication of Schach's work, a powerful theory
has been developed which gives conditions for weak convergence in such circumstances; see Stroock and Varadhan
\cite{sv}, Chapter 11. Like Schach's results, Theorem \ref{thm:sbm} below --- concerning the weak convergence of
projections of SBM --- also follows as a consequence of this general theory. (Despite this, we are not aware that the
result has been made known explicitly.) On the other hand, Theorem \ref{thm:lnnrw} --- concerning random walk on the
hypercube --- can not be deduced from the same theory. This is quite simply because an arbitrary projection of the
random walk on the hypercube is non-Markov for finite $d$. Hence, it is the proof of Theorem \ref{thm:lnnrw} that
occupies the better part of the sequel.





\section{Set-up and main results}

Unless otherwise stated, we continue to adopt the notational convention that vectors appear in bold typeface and the
value $x_j$ is assumed to be the $j$th component of a vector $\x$. In addition, a $\mc D$ above a binary relation
indicates that the relation holds in the sense of probability law, while $:=$ indicates a notational definition. Also,
$\N_0 := \N \cup \{0\}$ denotes the set of non-negative integers and for any $m\in\N$, $[m]:= \{1,\dots,m\}$. Any
convergence statements made in the sequel are intended to be understood with respect to the limit $d \to\infty$.

Define the OU-process $U = \{U_t\}_{t\ge 0}$ as the diffusion process with drift and diffusion coefficients given by
$b(u) := - u$ and $a(u) := 2$. That is, the infinitesimal generator of $U$ is given by \beqn \label{L} \mc L :=
\frac{\pdv^2 }{\pdv u^2} - u \frac{\pdv }{\pdv u}. \neqn Here and subsequently, we shall assume $U_0 = u$ is
deterministic so that $U$ is a Gaussian process.

Let $\X\sd = \{\X\sd_t\}_{t \ge 0}$ denote a random walk on $\Sd$ with $\X\sd_t$ representing the location of the
walker at time $t$. For each $d\in\N$ we choose a `direction' $\bth\sd \in \So$ and define the `projected process'
$Y\sd = \{Y\sd_t\}_{t\ge0}$ by \beqns Y\sd_t := \langle \bth\sd, \X\sd_t\rangle, \neqns where
$\langle\cdot,\cdot\rangle$ denotes the conventional inner product. Note that $Y\sd$ depends on $\bth\sd$ but that this
is not explicitly highlighted in the notation. We will assume that $\X\sd$ starts from a given initial position
$\X\sd_0 = \x\sd$ so that $Y\sd_0 = y\sd := \langle \bth\sd, \x\sd\rangle$ is deterministic.

Though we make the assumption that $\X\sd$ starts at the deterministic point $\X\sd_0 = \x\sd$, it will become evident
that analogous results hold when $\X\sd$ is a stationary random walk, {\it i.e.} $\X\sd_0$ is distributed uniformly on
the state-space. Commuting the role of randomness in $\X\sd$ and $\bth\sd$ implies a randomized central limit theorem
for the case when $\X\sd_0$ is again fixed but $\bth\sd$ is chosen uniformly from the state-space (and then
normalised). Note also that while our theorems will be stated only for $1$-dimensional projections of $\X\sd$, the
results are readily extended, via the Cram\'er-Wold device, to $K$-dimensional projections of $\X\sd$, $K<\infty$.


\subsection{Continuous case: spherical Brownian motion}

Here we take $\X^{(d)}$ to be SBM on $\Sd$. Using the definition stated in It\'o and McKean \cite{im}, SBM on $\Sd$ is
the unique diffusion process with infinitesimal generator:
\beqn \label{Delta} \Delta_d = \sum_{i \in [d]}\frac{\pdv^2}{\pdv x_i^2} - \frac 1{d}\sum_{i,j \in [d]} x_ix_j
\frac{\pdv^2}{\pdv x_i\pdv x_j} - \frac{d-1}{d} \sum_{i\in [d]}x_i\frac{\pdv}{\pdv x_i}. \neqn Alternative yet
equivalent characterisations of SBM are given in Stroock \cite{stroock} and Rogers and Williams \cite{rw}.

\begin{thm} \label{thm:sbm}
Let $\X\sd$ be SBM on $\Sd$. If $y^{(d)} \to u$, then \beqns Y^{(d)} \td U. \neqns
\end{thm}

\begin{proof}
By the symmetry of SBM it follows that $Y\sd$ is Markov. What's more, its infinitesimal generator is determined from
that of $\X^{(d)}$ simply by studying the action of $\Delta_d$ on functions dependent only on $y = \bth^{(d)}\cdot \x$.
From the expression given in \eqref{Delta}, it follows immediately that the infinitesimal generator of $Y^{(d)}$ is
given by \beqn \label{Ld} \mc L_d = -\frac{d-1}{d}y\frac{\pdv}{\pdv y} + \lrb 1 - \frac{y^2}{d} \rrb \frac{\pdv^2}{\pdv
y^2}. \neqn We read off the drift and diffusion coefficients as $b_d(y):= -(d-1)y/d$ and $a_d(y):= 2(1-y^2/d)$
respectively. Now, since ($i$): $U$ is the unique process started at $U_0 = u$ with infinitesimal generator $\mc L$;
($ii$): $a_d$ and $b_d$ are continuous and bounded uniformly in $d$ on compact subsets of $\R$; ($iii$): $a_d$ (resp.
$b_d$) converges pointwise to $a$ (resp. $b$) on compact subsets of $\R$; the result follows by Theorem 11.1.4, page
264, of Stroock and Varadhan \cite{sv}.

\end{proof}


\subsection{Discrete case: random walk on a hypercube}

From here on we take $\X\sd$ to be a simple, `lazy', nearest neighbour random walk (LNNRW) on the vertices of the
hypercube $\mc B^d$, where $\mc B := \{-1,1\}$. Two vertices in $\mc B^d$ are nearest neighbours if they differ in
exactly one coordinate. We assume that at regular clock pulses, separated by time intervals of length $\de :=
\frac{2p}d$, a LNNRW on $\mc B^d$ is `lazy' (remains stationary) with probability $1-p$, or moves to any given one of
its $d$ nearest neighbour vertices with equal probability $\frac pd$. We assume $p \in (0,1]$ may depend on $d$ (\eg $p
= d/(d+1)$).

\begin{thm} \label{thm:lnnrw}
Let $\X\sd$ be LNNRW on $\mc B^d$. If $y^{(d)} \to u$ and $|\bth\sd|_\infty \to 0$, then \beqn Y^{(d)} \td U. \neqn
\end{thm}

\begin{proof}
We apply the general program of Billingsley \cite{billingsley1}. The convergence of the finite dimensional
distributions of $Y\sd$ to those of $U$ is given by Lemma \ref{fddy}; tightness of the sequence $\{Y\sd\}_{d\in\N}$ is
established by Lemma \ref{lem:tightness}.
\end{proof}

\begin{rem}
The extra condition appearing in Theorem \ref{thm:lnnrw} that was absent from Theorem \ref{thm:sbm} is due to the lack
of complete spherical symmetry of the LNNRW. As an example of why some condition on $\bth^{(d)}$ is necessary, consider
the choice $\bth^{(d)} = (1,0,\dots,0) \in \So$ for each $d=1,2,\dots$. In this case it is clear that $Y^{(d)}$ is a
two-valued process and thus in no way can approach a diffusion limit as $d\ra\infty$.
\end{rem}

\begin{rem}
It is possible to generalise the LNNRW model along the lines of Schach's multivariate urn model without changing the
conclusions of Theorem \ref{thm:lnnrw}; see Remark \ref{rem:paramgen}.
\end{rem}



\section{Results for random walk on the cube}
\rnc{\Z}{\bs Z}

We begin with a concrete characterisation of LNNRW on $\mc B^d$. Since $\X^{(d)}_t$ gives the location of the random
walker at the (continuous) real time $t$, we will also adopt the alternative notation $\X^{(d)}(n) \equiv \X^{(d)}_t$,
$n\de \le t < (n+1)\de$, so that $\X^{(d)}(n)$ represents the location of the walker after $n$ clock pulses.
Here and subsequently let $X_1,X_2,\dots$ be a sequence of i.i.d. copies of $\X^{(1)}$ (LNNRW on the $1$-cube) and let
$\M^{(d)}(n) \sim \mult (n;\sfrac 1{d},\dots, \sfrac 1{d})$ denote a multinomial random vector with parameters
$(n;\sfrac 1{d},\dots, \sfrac 1{d})$. Since at any clock pulse the two choices of where to walk and whether to walk are
interchangeable, a moments's reflection will confirm that \beqn \label{xdist} \X\sd(n) \eqd \big(
X_1(M\sd_1(n)),\dots,X_d(M\sd_d(n)) \big). \neqn

Now, let $\Z\sd$ be a discrete-time random process on $\mc B^{d}$ with i.i.d. coordinate processes, each of which is
equal in distribution to any of the identically distributed, but dependent, coordinate processes of $\X\sd$. That is,
for each $n \in \N_0$, \beqn \label{zdist} \Z\sd(n) :\eqd \big (X_1(B\sd_1(n)),\dots,X_d(B\sd_d(n))\big ), \neqn where
$\bs B\sd(n)$ is a vector of independent $\Bi (n;\frac 1d)$ binomial random variables each with parameters $(n;\frac
1d)$. The proximity of the moments of the finite dimensional distributions of $\X\sd$ to those of $\Z\sd$ will be the
result that's useful in the sequel. Now we introduce some notation that helps us to precisely state and prove the
result we require.

Due to the symmetry of $\mc B^d$ and the arbitrariness of $\bth\sd$, we may, without loss of generality, restrict our
attention to only a single choice of initial position. Thus, we henceforth assume that $\X\sd_0$ starts at $\x\sd =
(1,\dots,1) \in \mc B^d$.

Pre-empting the treatment of the finite dimensional distributions of $Y\sd$, let $0 = t_0 < t_1 < \dots < t_K < \infty$
be a sequence of $[0,\infty)$-valued times and let $n_0,n_1,\dots,n_K$ be the corresponding sequence of $\N_0$-valued
`$\de$-counts' such that $n_k$ is the integer part of $t_k/\de$; note that $n_k$ depends on $d$. From here on we shall
often refrain from indicating the dependence on $d$ explicitly with the superscript $\sd$. We shall also utilize the
shorthand $\bs V(n_k) \equiv \bs V_k$ for any vector-valued process $\bs V \in \R^d$.

For each $k\in [K]$, let $\M_k' \sim \mult (n_k-n_{k-1};\frac 1d,\dots, \frac 1d)$ and $\bs B_k'$ be a vector of i.i.d.
$\Bi (n_k-n_{k-1};\frac 1d)$ random variables, such that $\M_k = \sum_{j=1}^k \M_j'$ and $\bs B_k = \sum_{j=1}^k \bs
B_j'$, and introduce \beqns \X_k' :\eqd ( X_1(M_{k1}'),\dots,X_d(M_{kd}')), \quad \Z_k' :\eqd
(X_1(B_{k1}'),\dots,X_d(B_{kd}')).\neqns In view of the Markov property we may deduce that for any $i$: \beqns
X_{ki}\eqd \prod_{j=1}^k X_{ji}', \qquad Z_{ki} \eqd \prod_{j=1}^k Z_{ji}', \neqns assuming that the $\X_k'$ (resp.
$\Z_k'$) are mutually independent.

Now, fix a multi-index $i = (i_1,\dots,i_L) \in [d]^L$ and constants $l_1,\dots,l_K$ such that $l_1+\dots + l_K =
L$. For each $k \in [K]$, let $L_k := l_1+\dots +l_k$, $L_k' := L_{k-1}+1$ and $J_k(i)$ be the set containing precisely those $j \in [d]$
that occur with odd multiplicity in the multi-index $(i_{l_k'},\dots,i_L)$ of length $l_k+\dots +l_K$. We will also
need $\eta_k(i) := |J_k(i)|$.

\begin{lem} \label{lem:XZmoms}
For any multi-index $i = (i_1,\dots,i_L) \in [d]^L$,
\begin{align}
\label{eq:EprodX} \Exp \big \{ \prod_{k=1}^K \prod_{l=L_k'}^{L_k} X_{ki_l} \big \} &=
\prod_{k=1}^K ( 1 - \eta_k(i) \de)^{n_k-n_{k-1}},\\
\label{eq:EprodZ} \Exp \big \{ \prod_{k=1}^K \prod_{l=L_k'}^{L_k} Z_{ki_l} \big \} &= \prod_{k=1}^K ( 1 -
\de)^{\eta_k(i)(n_k-n_{k-1})}.
\end{align}
\end{lem}

\begin{proof}
We prove only \eqref{eq:EprodX}, the proof of \eqref{eq:EprodZ} is analogous. Fix $i$ and set $J_k \equiv
J_k(i)$.\beqns \Exp \big \{ \prod_{k=1}^K \prod_{l=L_k'}^{L_k} X_{ki_l} \big \}= \Exp \big \{ \prod_{k=1}^K
\prod_{l=L_k'}^{L_k} \prod_{j=1}^k X_{ji_l}' \big\} = \prod_{j=1}^K \Exp \big \{ \prod_{k=j}^K \prod_{l=L_k'}^{L_k}
X_{ji_l}'\big \} = \prod_{k=1}^K \Exp \big \{ \prod_{l=L_k'}^{L} X_{ki_l}'\big \}.  \neqns Now, using the fact that
$(X_{ki}')^r = X_{ki}'$ if $r$ is odd and $(X_{ki}')^r =1$ otherwise, we see that \beqn \label{eq:EprodXa} \Exp \big \{
\prod_{l=L_k'}^{L} X_{ki_l}'\big \} = \Exp \big\{ \prod_{j \in J_k} \Exp \big\{ \X^{(1)}(M_{kj}') \mid M_{kj}' \big\}
\big\} = \Exp \big\{ \prod_{j \in J_k} \la^{M_{kj}'} \big\}, \neqn where $\la := 1-2p$ is the non-unit eigenvalue of
the transition probability matrix of $\X^{(1)}$. The result now follows by noting that $\Exp \big\{ \prod_{j \in J_k}
\la^{M_{kj}'} \big\}$ is the probability generating function, evaluated at $\la$, of $\sum_{j \in J_k} M_{kj}' \sim \Bi
\lrb n_k - n_{k-1}; \frac{\eta_k(i)}d \rrb$.
\end{proof}

\begin{rem}\label{rem:paramgen}
All quantitative information that is used in the proof of Theorem \ref{thm:lnnrw} can be traced back to Lemma
\ref{lem:XZmoms}, which itself hinges on certain independence assumptions and the valuation of $\Exp \big\{ \prod_{j
\in J_k} \la^{M_{kj}'} \big\}$ in \eqref{eq:EprodXa}. Thus, generalisations of the LNNRW model that leave Theorem
\ref{thm:lnnrw} unchanged become apparent. For example, with reference to \eqref{xdist}, if we exchange $\M\sd(n)$ for
$\bs N\sd (n) \sim \mult (n;\bs\phi\sd)$ and the i.i.d. $\mc B$-valued LNNRW's $X_1,X_2,\dots$ for the independent $\mc
B$-valued recurrent Markov chains $V\sd_1,V\sd_2,\dots$, the respective transition probability matrices of which have
non-unit eigenvalues $\la\sd_1,\la\sd_2,\dots$, then provided there exists a $\de$ dependent on $d$ such that
$\phi\sd_j(1-\la\sd_j)\de^{-1} \to 1$ for each $j$, the conclusion of Theorem \ref{thm:lnnrw} remains valid.
\end{rem}


\subsection{Convergence of finite dimensional distributions} \label{convYcube}

\begin{lem} \label{fddy}
If $|\bth\sd|_\infty \to 0$ and $y\sd \to u$, then \beqns ( Y\sd_{t_1},\dots,Y\sd_{t_K}) \td
(U_{t_1},\dots,U_{t_K}).\neqns
\end{lem}

\begin{proof}

By the Cram\'er-Wold device (Billingsley \cite{billingsley1}, Theorem 7.7) it is enough to show that for any $\bs \phi \in \R^K$, \beqns \Psi_d := \sum_{k=1}^K \phi_k Y\sd_{t_k} = \sum_{k=1}^K \phi_k \langle
\bth\sd, \X\sd_k\rangle \td \Ga := \sum_{k=1}^K \phi_k U_{t_k}. \neqns Since the (Gaussian) law of $\Ga$ is uniquely
determined by its sequence of moments, this can be achieved through a method of moments argument (Gut \cite{gut05}, p.
237) by showing that, \beqn \label{eq:mmoms} \Exp \lcb \Psi_d^L \rcb \to \Exp \lcb \Ga^L \rcb, \qquad \; L = 1,2,\dots.
\neqn

Now, let $\Upsilon_d := \sum_{k=1}^K \phi_k \langle \bth\sd,\Z\sd_k\rangle = \sum_{j \in [d]}\tht\sd_j \xi\sd_j$ be the sum of $d$ independent random variables, with $\xi\sd_j := \sum_{k=1}^K\phi_k Z\sd_{kj}$, and set $\si_d^2 := \Var \{\Upsilon_d\}$.
By taking appropriate linear combinations of formula \eqref{eq:EprodZ} with $K=1,2$ and $L=1,2$, it is straightforward to show that the first two moments of $\Upsilon_d$ converge to those of $\Ga$. Moreover, the fact that $|\bth\sd|_2 = 1$ and the condition that $|\bth\sd|_\infty \to 0$ is enough to ensure there exists an $r>2$ such that \beqns \sum_{j=1}^d \si_d^{-r}\Exp |\tht\sd_j\xi\sd_j-\Exp\{\tht\sd_j\xi\sd_j\} |^{r} = \si_d^{-r}\Exp |\xi\sd_1-\Exp\{\xi\sd_1\} |^{r} \sum_{j=1}^d |\tht\sd_j|^r \ra 0, \neqns so that the Lyapounov condition (Gut \cite{gut05}, p. 339) is satisfied. Hence we may conclude from Lyapounov's central limit theorem that $\Upsilon_d \to \Ga$. What's more, we may use the Marcinkiewicz-Zygmund inequalities (Gut \cite{gut05}, p. 146) to verify that $\Exp \{\Upsilon_d^L\}$
is bounded for each $L\in\N$, implying $\Upsilon_d^L$ is uniformly integrable for each $L \in\N$ and thus $\Exp \{\Upsilon_d^L\} \to \Exp \{\Ga^L\}$ for each $L\in \N$ (Billingsley
\cite{billingsley1}, Theorem 5.4). Finally, since $\Exp \lcb \Psi_d^L - \Upsilon_d^L \rcb \to 0$ for each $L\in\N$
(Lemma \ref{lem:uptopsi} below), we conclude that $\Exp \{\Psi_d^L\} \to \Exp\{ \Ga^L\}$, as required.
\end{proof}

Before coming to the proof of Lemma \ref{lem:uptopsi} cited above, we need to cover an intermediary result. Recall that
a set $\pi$ of non-empty subsets of a finite set $\Sigma$ is a partition of $\Sigma$ if the elements of $\pi$ are mutually
disjoint {\it and} $\Sigma = \cup_{\si \in \pi} \si$. For two partitions $\pi = \{\pi_1,\dots,\pi_m\}$ and $\nu =
\{\nu_1,\dots,\nu_l\}$ of the same finite set, we will write $\nu \prec \pi$ if $l < m$ {\it and} each $\nu_j$ is a
union of $\pi_j$'s. We write $\nu \preceq \pi$ if either $\nu \prec \pi$ or $\nu = \pi$.

To every multi-index $i = (i_1,\dots,i_L) \in [d]^L$, or equivalently, mapping $i:[L] \to [d]:l \mapsto i_l$, there
corresponds a partition of $[L]$: \beqns \nu_i := \{i^{-1}(k) \subset [L]: k \in [d]\} \setminus \{\emptyset\}, \neqns
where $\emptyset$ denotes the empty set and $i^{-1}(k) := \{l \in [L]:i_l = k\}$ is the pre-image of $k$ under the
mapping $i$. For any partition $\pi = \{\pi_1,\dots, \pi_m\}$ of $[L]$, define \beqns I_{\pi} := \{i \in [d]^L: \nu_i =
\pi\}, \neqns so that for each $i \in I_{\pi}$ there exists $m$ {\it distinct} numbers $b_1,\dots,b_m$ such that for
each $k \in [m]$, $i_l = b_k$ for all $l \in \pi_k$. In addition, define $I_{\prec \pi} = \cup_{\nu \prec \pi} I_{\nu}$
and $I_{\preceq \pi} = I_{\pi} \cup I_{\prec \pi}$ so that, in particular, $I_{\preceq \pi}$ contains precisely those
$i \in [d]^L$ where there exists {\it not-necessarily distinct} numbers $b_1,\dots,b_{m}\in [d]$ such that for each $k
\in [m]$, $i_l = b_k$ for all $l \in \pi_k$.


\begin{prop} \label{prop:dtoLsum}
Given an array $(a_{jl}): j \in [d], l \in [L]$, set $A_s := \sum_{j \in [d]}\prod_{l \in s} a_{jl}$ for any $s \subset
[L]$. Let $(c_{\pi,\nu})_{\nu \preceq \pi}$ be the triangular array of constants, indexed by partitions of $[L]$, that
satisfies the recursion: $c_{\pi,\nu} = - \sum_{\nu \preceq \mu \prec \pi}c_{\mu,\nu}$, for $\nu \prec \pi$, and $c_{\pi,\pi} = 1$. Then, \beqn
\label{eq:A21} \sum_{i \in I_{\pi}} \prod_{l \in [L]} a_{i_ll} = \sum_{\nu \preceq \pi} c_{\pi,\nu} \prod_{s \in
\nu}A_s. \neqn
\end{prop}

\begin{proof}
The proof is via induction on $|\pi|$. As the first step: when $|\pi|=1$ we must have $\pi = \{[L]\}$ and $I_{\pi} =
\{i \in [d]^L: \exists k \in [d] \mbox{ s.t. } i = (k,\dots,k)\}$ so that \beqns \sum_{i \in I_{\pi}} \prod_{l \in
[L]}a_{i_ll} = \sum_{k \in [d]} \prod_{l \in [L]}a_{kl} = A_{[L]}. \neqns

Assuming the induction hypothesis (\ref{eq:A21}) to be true for any partition $\pi$ of $[L]$ such that $|\pi| \le m$,
we now proceed with the induction step. If $\pi = \{\pi_1,\dots,\pi_{m+1}\}$ is a partition of $[L]$, then
\begin{align*}
\sum_{i \in I_{\pi}} \prod_{l\in [L]} a_{i_ll} &= \sum_{i \in I_{\preceq \pi}} \prod_{l \in [L]}a_{i_ll} - \sum_{i
\in I_{\prec \pi}} \prod_{l\in [L]}a_{i_ll}\\
&= \sum_{i \in [d]^{m+1}} \prod_{k=1}^{m+1} \prod_{l \in \pi_k} a_{i_kl} - \sum_{\nu \prec \pi} \sum_{i \in I_{\nu}}
\prod_{l\in [L]}a_{i_ll}\\
&= \prod_{k=1}^{m+1}\sum_{j \in [d]} \prod_{l \in \pi_k} a_{jl} - \sum_{\nu \prec \pi} \sum_{\mu \preceq \nu} c_{\nu,\mu} \prod_{s \in \mu} A_s\\
&= \prod_{k=1}^{m+1} A_{\pi_k} -  \sum_{\mu \prec \pi} \big ( \sum_{\mu \preceq \nu \prec \pi} c_{\nu,\mu} \big
)\prod_{s \in \mu} A_s\\
&=  \sum_{\mu \preceq \pi} c_{\pi,\mu}' \prod_{s \in \mu} A_s,
\end{align*}
where $c_{\pi,\pi}' = 1$ and $c_{\pi,\mu}' = - \sum_{\mu \preceq \nu \prec \pi}c_{\nu,\mu}$. Hence $c_{\pi,\nu}' =
c_{\pi,\nu}$ as required.
\end{proof}

\begin{rem} \label{rem:A21}
Note that the number of terms in the sum on the left hand side of (\ref{eq:A21}) depends only on $d$ whereas the number
of terms in the sum on the right hand side of the same equation depends only on $L$.
\end{rem}

\begin{lem}\label{lem:uptopsi}
For each $L \in \N$, $\Exp \lcb \Psi_d^L - \Upsilon_d^L \rcb \to 0$.
\end{lem}

\begin{proof}
Let $\bs l \in \N_0^K$ be such that $l_1 + \dots + l_K = L$. Define the shorthand \beqns E_d(\bs l) := \Exp \blcb
\langle\bth, \X_1\rangle^{l_1}\dots \langle\bth, \X_K\rangle^{l_K}\brcb - \Exp \blcb \langle\bth,
\Z_1\rangle^{l_1}\dots \langle\bth, \Z_K\rangle^{l_K}\brcb \neqns and introduce $\Pi_L := \{I_{\pi}: \pi \mbox{ is a
partition of $[L]$}\}$ as a partition of $[d]^L$ into disjoint subsets of multi-indices over which $\bs{\eta}(i)$
remains constant. Agree to allow $\bs{\eta}(\pi) \equiv \bs{\eta}(i)$ whenever $i \in I_{\pi}$. Now, after using a
multinomial expansion we get\beqn |\Exp \lcb \Psi_d^L - \Upsilon_d^L\rcb|
\le \lrb |\phi_1|+\dots + |\phi_K|\rrb^L \max_{\bs l} |E_d(\bs l)|, \label{eq:upsL}\neqn so that it suffices to show
that $E_d(\bs l) \to 0$ for each admissable $\bs l$.

Using Lemma \ref{lem:XZmoms}, we have
\begin{align}
E_d(\bs l) &= \sum_{i \in [d]^L} \Big [ \Exp \blcb \prod_{k=1}^K \prod_{l = L_k'}^{L_k} \tht_{i_l}X_{ki_l} \brcb -
\Exp \blcb \prod_{k=1}^K \prod_{l = L_k'}^{L_k} \tht_{i_l}Z_{ki_l} \brcb \Big ]\nn\\
&= \sum_{\pi \in \Pi_L} f_d(\pi) \blrb \sum_{i \in I_{\pi}} \prod_{l=1}^L \tht_{i_l} \brrb, \label{eq:Edlast}
\end{align}
where for all $\pi \in \Pi_L$, \beqns f_d(\pi) := \prod_{k=1}^K (1-\eta_k(\pi)\de)^{n_k - n_{k-1}} - \prod_{k=1}^K
(1-\de)^{\eta_k(\pi)(n_k - n_{k-1})} \to 0. \neqns By assumption: $\sum_{j \in [d]} \tht_j \to u$, $\sum_{j \in [d]}
\tht_j^2 =1$ and for all $r>2$, $\big |\sum_{j \in [d]} \tht_j^r\big | \le \max_{j \in [d]}|\tht_j|^{r-2}\sum_{j\in[d]}|\tht_j|^2
\to 0$; hence, replacing $a_{jl}$ with $\tht_j$ for all $l$, we conclude via Proposition \ref{prop:dtoLsum} that
$\sum_{i \in I_{\pi}} \prod_{l=1}^L \tht_{i_l}$ remains bounded as $d \to \infty$; see Remark \ref{rem:A21}. Noting
that $|\Pi_L|$ is independent of $d$, it follows that the sum (\ref{eq:Edlast}) consists of a fixed number of terms
each of which tends to zero as $d \to \infty$, thus completing the proof.
\end{proof}


\subsection{Tightness}

\begin{lem} \label{lem:tightness}
If $y\sd \to u$, then $\{Y\sd\}_{d = 1}^\infty$ is tight.
\end{lem}

\begin{proof}
We appeal to the tightness criterion of Theorem 4.1, page 355, Jacod and Shiryaev \cite{js03}. Let $\vep
> 0$ and let $C$ denote a generic constant independent of $d$ and $t_3$. Since there is zero probability of a jump discontinuity at time zero,
it suffices to show that \beqn \label{tightdeqn} \Prb \lrb |Y\sd_{t_3}-Y\sd_{t_2}| \ge \varepsilon,
|Y\sd_{t_2}-Y\sd_{t_1}| \ge \vep \rrb \le \frac{C}{\vep^{3}} (t_3-t_1)^{\frac 32}; \qquad d=2,3,\dots . \neqn To
establish this result we follow the same method as used on page 459 of Schach \cite{schach}. In what follows, we assume
$d\ge 2$ and drop the $\sd$ notation once again. 

\begin{lem} \label{lem1} 
\beqn \Prb \lrb |Y_{t_2}-Y_{t_1}| \ge \varepsilon \; | Y_{t_1} \rrb \le \frac{4}{\varepsilon^2} (n_2-n_1)\de \lrb 1 +
Y_{t_1}^2 \rrb. \neqn
\end{lem}

\begin{proof}
Define $\nu_k(x)= 1-kx\de$ and let \beqns f(x) = 1 - \nu_2^{n_2-n_1}(x) + Y_{t_1}^2\lcb \nu_2^{n_2-n_1}(x) -
2\nu_1^{n_2-n_1}(x) + 1 \rcb. \neqns Then by the Mean Value Theorem, \beqn f(1)-f(0) = \Exp \lcb \lrb Y_{t_2}-
Y_{t_1}\rrb^2 | Y_{t_1} \rcb \le 2(n_2-n_1)\de \lsb 1 + 2Y_{t_1}^2\rsb. \neqn The result now follows by Chebyshev's
inequality.
\end{proof}

\begin{lem}\label{lem2}
If $y\sd \le C$, then for any $L\in \N$, \beqns |\Exp \{Y_{t_3}^L\}|  \le C.\neqns
\end{lem}

\begin{proof}
First, if $d \le L$ we may use the trivial bound $|\Exp \{Y_{t_3}^L\}|  \le L^{\frac L2} \le C$. If, on the other hand, $d > L$, then for any $i \in [d]^L$, $0 \le \eta_3(i)\de \le 2$. Thus, following analogous arguments to those used in the proof
of Lemma \ref{lem:uptopsi}, we see that \beqns |\Exp \{Y_{t_3}^L\}| = \big |\sum_{i \in [d]^L} \blrb
\prod_{l=1}^L \tht_{i_l} \brrb (1-\eta_3(i)\de)^{(n_3 - n_{2})}\big | \le \sum_{\pi \in \Pi_L} \big | \sum_{i \in
I_{\pi}} \prod_{l =1}^L \tht_{i_l} \big | \le C.\neqns 
\end{proof}

\begin{lem}
\beqn \label{whynot} \Prb \lrb |Y_{t_3}-Y_{t_2}| \ge \varepsilon, |Y_{t_2}-Y_{t_1}| \ge \vep \rrb \le
\frac{C}{\vep^{3}}(n_3-n_2)(n_2-n_1)^{\frac 12}\de^{\frac 32}. \neqn
\end{lem}


\begin{proof} Utilising the Markov property of $\X$, the
Cauchy-Schwarz inequality and Lemmas \ref{lem1} and \ref{lem2}, we obtain:
\begin{align*}
\Prb \{ |Y_{t_3}-Y_{t_2}| \ge &\vep, |Y_{t_2}-Y_{t_1}| \ge \vep \}\\[3pt]
&= \Exp \lsb \Prb \lcb |Y_{t_3}-Y_{t_2}| \ge \varepsilon, \;|Y_{t_2}-Y_{t_1}| \ge \vep \;|\X_{t_2}\rcb \rsb \\[3pt]
&= \Exp \lsb \Prb \lcb |Y_{t_3}-Y_{t_2}| \ge \varepsilon \;|\X_{t_2} \rcb \Prb \lcb |Y_{t_2}-Y_{t_1}| \ge \vep\; |\X_{t_2}\rcb \rsb \\[3pt]
&\le \frac{C(n_3-n_2)\de}{\vep^2}\Exp \lsb \Prb \lcb |Y_{t_2}-Y_{t_1}| \ge \varepsilon \;|\X_{t_2} \rcb
\lrb 1 +  Y_{t_2}^2\rrb \rsb  \\
&\le \frac{C(n_3-n_2)\de}{\vep^2}\lsb \Exp \lsb \Prb \lcb |Y_{t_2}-Y_{t_1}| \ge \varepsilon |\X_{t_2} \rcb\rsb^2
\Exp\lsb 1 +  Y_{t_2}^2\rsb^2\rsb^{\frac 12}  \\
&\le \frac{C(n_3-n_2)\de}{\vep^2} \lsb \Prb \lcb |Y_{t_2}-Y_{t_1}| \ge \varepsilon \rcb \rsb^{\frac 12} \\
&= \frac{C(n_3-n_2)\de}{\vep^2} \lsb \Exp\lsb\Prb \lcb |Y_{t_2}-Y_{t_1}| \ge \varepsilon \;|Y_{t_1} \rcb \rsb\rsb^{\frac 12} \\
&\le \frac{C(n_3-n_2)(n_2-n_1)^{\frac 12}\de^{\frac 32}}{\vep^{3}}
\lsb \Exp\lsb 1 +   Y_{t_1}^2\rsb\rsb^{\frac 12}  \\
&\le \frac{C(n_3-n_2)(n_2-n_1)^{\frac 12}\de^{\frac 32}}{\vep^{3}}.
\end{align*}
\end{proof}


We now verify that \eqref{whynot} implies the tightness condition \eqref{tightdeqn}. Suppose first that $n_3>n_2>n_1\ge
0$. Then clearly $n_3-n_2\ge 1$, $n_2-n_1\ge 1$ and $n_3-n_1 \ge 2$ from which it follows immediately that \beqns
\label{final2} (n_3-n_2)(n_2-n_1)^{\frac 12} \le (n_3-n_1-1)^{\frac 32} \le (t_3-t_1)^{\frac 32} \de^{-\frac 32}.
\neqns Moreover, if $n_3=n_2$ and/or $n_2=n_1$, the above inequality is trivially satisfied and thus it in fact holds
for $n_3\ge n_2\ge n_1 \ge 0$.

\end{proof}


\section*{Acknowledgements}

The author would like to thank Terry Lyons for suggesting the problem and Stephen Buckley, Svante Janson and Gesine
Reinert for helpful suggestions that led to improvements of the proofs.




\bibliography{phdrefs}

\end{document}